\documentclass[12pt,a4paper]{article}

\usepackage{amsmath,amssymb,amsfonts,amsthm,latexsym,amscd}
\usepackage[english]{babel}
\usepackage{graphics}
\usepackage{wasysym}
\usepackage{framed}
\usepackage[mathscr]{eucal}   	
\usepackage{color}
\usepackage[all,2cell]{xy} \UseAllTwocells 		% for diagrams + 2-cells 	
\definecolor{darkjunglegreen}{rgb}{0.1, 0.14, 0.13}
\definecolor{mydarkgreen}{rgb}{0.0, 1.0, 0.2}

\newcommand{\C}{\mathcal{C}}
\newcommand{\D}{\mathcal{D}}
\newcommand{\E}{\mathcal{E}}\newcommand{\F}{\mathcal{F}}\newcommand{\G}{\mathcal{G}}
\renewcommand{\H}{\mathcal{H}} 	%re!
 
\renewcommand{\L}{\mathcal{L}}		%re!

\newcommand{\Df}{{\scr{D}}}

\DeclareSymbolFont{rsfs}{U}{rsfs}{m}{n}
\DeclareSymbolFontAlphabet{\scr}{rsfs}% was ... mathscr

\DeclareMathOperator{\Ob}{Ob}

\DeclareMathOperator{\id}{id}

\newtheorem{theorem}{Theorem}[section]				% theorems numbered under sections

\theoremstyle{definition}        
\newtheorem{definition}[theorem]{Definition}
\newtheorem{remark}[theorem]{Remark}

\title{\textbf{\LARGE{Two partial monoid structures on a set}}}
\author{
\normalsize Rachel A.D. Martins 
\date{2 February 2015}
\\
\normalsize \textit{Centro de \'Algebra da Universidade de Lisboa,}
\\
\normalsize \textit{Av. Prof. Gama Pinto, 2,}
\\
\normalsize \textit{1649-003 Lisboa, Portugal}
}

\begin{document}

\maketitle

\begin{abstract}
It is well-known that small categories have equivalent descriptions as partial monoids. We
provide a formulation of partial monoid and partial monoid homomorphism involving $s$ and $t$
instead of identities and then following
a recent investigation into involutive double categories, we prove that a double
category is equivalent to a set equipped with two partial monoid structures in which all structure maps
$(s,t,\circ)$ are partial monoid homomorphisms. We discuss the light this purely algebraic perspective
sheds on the symmetries of these structures and its applications.
Iteration of the procedure leads also to a pure algebraic formulation of the notion of $n$-fold
category.
\end{abstract}

\section{Introduction}

The partial monoid description of a small category \cite{ML} can be thought of as another
generalisation of the notion of group. The key idea in translating higher categorical structures into a
pure algebraic language, is to afford the $m < n$
cells the status of $n$-cells, so in a purely algebraic double category or 2-category, the objects
and 1-morphisms are described as members of the set of 2-morphisms. The author learned
about this point of view from discussions with Pedro Resende, Paolo Bertozzini and Roberto Conti and
the original theme dates back to some handwritten notes by John Roberts in the 1970s.

\medskip

Let $\C$ be a category. In this contribution following \cite{idc}, we replace
the
identity map $\id_{\C} : \Ob(\C) \to \C$ with source and target maps $s,t : \C \to \C$
leading to an \emph{emergence} of $\Ob(\C)$ as a (not pre-distinguished) subset of the set of morphisms
$\C$. In our construction, the notion of identity morphism does not appear anywhere explicitly because
it is completely absorbed into the data provided by $s$ and $t$. We provide a concise translation of the
notion of double category  (a
category internal in the category of small categories \cite{E1st}) into a pure
algebraic language, in terms of a set with
two commuting partial monoid structures. Note that Ehresmann's definition of double category is
very concise since it takes account of all the symmetries involved in the structure. The notion given
below of partial monoid homomorphism, which contains 3 conditions (the identity preservation condition
of a covariant functor is split into 2 conditions, one of which involves the source map
$s$ and the other, the  target map $t$) is the key factor that allows our translation to be almost as
concise.  The technique can
be iterated $n$ times, and leads to a corresponding purely algebraic formulation of $n$-fold category
\cite{BH}.

\medskip

As detailed in \cite{idc}, the notion of double
category
together with
its 4-fold system of opposite categorical structures, can be illustrated by any one choice
of 4 different classes of cubical set, (or square 2-cell), although only one of these choices exactly
corresponds to the original Ehresman notion of double category. The particular
construction considered here, results in the elimination of all the other three choices of cubical set
class. (The term ``double partial monoid'' has been reserved for the pure algebraic formulation of
double categories including the full generality of allowing all such choices. A detailed expanded
definition of double partial monoid will appear in a coming preprint \cite{idc}.) 

\medskip

In this language, the notion of 2-category has a more complicated translation than the notion of
double category. This underlines the point of view that double categories are more naturally occurring
algebraic objects.

\medskip

We introduce some of our motivations in this paragraph and discuss some applications in the
final section. (i) There are several theories relating inverse semigroups with topological groupoids
(including
\cite{Pedro},
\cite{Lawson 1}). Viewing groupoids as partial monoids with all inverses, might allow a comparison to
be
made within a single algebraic language. In contrast to working in Top to construct
groupoids, where one usually begins with a space $X$ and next a space of arrows is defined over $X$,
examples of localic groupoids \cite{Pedro} tend to arise from a set of arrows with the structure of a
groupoid. Consider for
example the groupoid of ultrafilters of a tiling semigroup \cite{Tiling smgs},
\cite{Kellendonk}, where the groupoid units are constructed directly and afterwards are interpreted as
objects. Other applications to algebra may follow from comparing congruences of higher
categories viewed as sub-structures of sets equipped with several commuting partial monoid structures.
(ii) Double groupoids have proved to be important tools in homotopy theory \cite{BH1} and
therefore it is worth looking for ways to generalise them so that some of the concepts can be explored
in a non-commutative context\footnote{For example, Prof. Ronald Brown has opened a popular discussion on
the construction of convolution C*-algebras for double groupoids, which has involved viewing a double
groupoid in terms of its underlying set of 2-arrows.}. (iii) In \cite{RR} Roberts and Ruzzi
repackage the notion of $n$-category in terms of a set of $n$-arrows subject to a list of axioms to
clarify
their approach to the non-Abelian cohomology of a poset.

%It is expected for example that difficult to visualize problems involving
%$n$-morphisms may become more tractable from such a viewpoint.
%One of the motivations for translating higher categorical
%structures into a pure algebraic language is so that their symmetries become more
%transparent.

\medskip

We will use the notation $(\C, \C^0, \iota, \circ, s,t)$ to refer to a category $\C$ where $\C^0$
denotes the objects $\Ob(\C)$ of $\C$, with $\id_{\C}$ or $\iota$ for the unit map $\iota : \C^0 \to
\C$ and $s$ and $t$ the source and target maps $s,t : \C \to \C^0$. For composable arrows $x$
and $y$ in $\C$, we write $(x,y) \in \C \times_{\circ} \C$.

\section{Partial monoids}

\begin{definition}\cite{idc}
A \textbf{\emph{partial monoid}} $(\C, \circ, s, t)$ is a set with a partially defined
associative\footnote{by associative in this context, we mean whenever one of the two terms $z \circ (x
\circ y)$ and $( z \circ x) \circ y$ exists,
the other one exists as well, and they coincide.}
multiplication
$\circ$
together
with two idempotent maps $s,t:\C \to \C$ called \emph{partial identities} such that 
$x \circ y$ is defined in $\C$ exactly when $s(x) = t(y)$ (and $z \circ x$ is defined in $\C$ exactly
when $t(x) = s(z)$),
% whenever 
%$s(x) = t(y)$, then $x \circ y$ is defined in $\C$, 
% whenever 
%$t(x) = s(z)$, then $z \circ x$ is defined in $\C$,
 and such that
$x \circ s(x) = x$, ~
$t(x) \circ x = x$,~
for all  $x \in \C$.
\end{definition}

\begin{remark}
The previous definition is equivalent to that of partial monoid by Mac Lane in \cite{ML}
because the category units (or identity morphisms) are obtained from the images of the source $s$ and
target $t$ partial identity maps whenever $s(x)=t(y)$.
\end{remark}

\emph{Caveat:} Given the partial monoid description of a category $\C$, it is not always
possible to fully reconstruct the objects $\Ob(\C)$ of the original category
$\C$. \footnote{If two definitions are equivalent, it does not follow that any
example can be precisely reconstructed from the other definition. This point is illustrated by recalling
that the real algebras $\mathbb{H}$ and $M_2(\mathbb{R})$ have the same unit, so in a category whose
objects are vector spaces there can be an ambiguity.}

\medskip

In this pure algebraic language, a covariant functor has an equivalent description as a partial monoid
homomorphism:

\begin{definition}
Let $(\C_1,\circ_1,s_1,t_1)$ and $(\C_2, \circ_2,s_2,t_2)$ be two partial monoids. A
\textbf{\emph{partial monoid homomorphism}} is a map $\F: \C_1 \to \C_2$ such that for each ~$x \circ_1
y ~ \in \C_1$,
\begin{gather}
 \F(x \circ_1 y) = \F(x) \circ_2 \F(y), \\
 \F(s_1(x)) = s_2(\F(x)),~~~~~~~~ \F(t_1(x)) = t_2(\F(x)).
\end{gather}
\end{definition}
The above contains 3 explicit conditions, whereas the definition of covariant functor contains only 2
lines.

\section{Double categories and partial monoids}

Recall that Ehresmann's construction of a \textbf{\emph{double category}}, is a category internal in the
category Cat of small categories \cite{E1st}. One transparent albeit simplified way to unpack this
definition, is to present a double category as a small category $\D$
with four categorical structures $(\D, \D^1_v, \circ^2_v, \iota^1_v, s^1_v,t^1_v)$,  $(\D,
\D^1_h, \circ^2_h, \iota^1_h, s^1_h,t^1_h)$, $(\D^1_v, \D^0, \circ^1_v, \iota^1_v, s^1_v,t^1_v)$ and
$(\D^1_h, \D^0, \circ^1_h, \iota^1_h, s^1_h,t^1_h)$, such that the partially defined associative
multiplications are compatible with all units $\iota$, and with one another (more details below). The
definition of double
category is expanded in various equivalent presentations in \cite{idc} and a good definition also
appears in \cite{BC}.

 \medskip
 
It is convenient to use cubical set diagrams \cite{L} to illustrate this.
According to \cite{idc}, there exist four equivalent classes of such square cell diagrams, each of which
captures a notion of double category, and each of which belongs
to a system of 4 mutually opposite categorical structures. Only one of these cubical set classes
represents the categorical compositions in an intuitive way, and this corresponds to the
original Ehresmann double category (recalled above).

\medskip

%\begin{enumerate}
% \item $s_h ( s_v) = s_v (s_h)$, ~ $t_h (t_v) = t_v (t_h)$, ~ $t_h (s_v) = s_v (t_h)(x)$, ~ $s_h
%(t_v) = t_v (s_h)$.
% \item
% \item
% \item
%\end{enumerate}

\noindent
\textbf{\emph{Class one:}} ~~~$s_h ( s_v) = s_v (s_h)$, ~ $t_h (t_v) = t_v (t_h)$, ~ $t_h (s_v) = s_v
(t_h)$, ~ $s_h (t_v) = t_v (s_h)$. \\
 \textbf{\emph{Class two:}} ~ $s_h ( s_v) = t_v (s_h)$, ~ $t_h (t_v) = s_v (t_h)$, ~ $t_h (s_v) = t_v
(t_h)$, ~ $s_v (s_h) = s_h (t_v)$.   \\
 \textbf{\emph{Class three:}} $s_h ( s_v) = s_v (t_h)$, ~ $t_h (t_v) = t_v (s_h)$, ~ $t_h (s_v) = s_v
(s_h)$, ~ $t_v (t_h) = s_h (t_v)$.\\
 \textbf{\emph{Class four:}} ~$s_h ( s_v) = t_v (t_h)$, ~ $t_h (t_v) = t_v (t_h)$, ~ $t_h (s_v) = t_v
(s_h)$, ~ $s_v (t_h) = s_h (t_v)$.
 
 \medskip

In the context of cubical set class one, the compatibility rules for the multiplications $\circ^2_h$ and
$\circ^2_v$
are expressed in symbols as follows:
\begin{itemize}
\item
 $(x \circ^2_v y) \circ^2_h (z \circ^2_v w) = (x \circ^2_h z) \circ^2_v (y \circ^2_h w)$~ whenever
both sides are defined,
 \item 
$s_v(x) \circ^2_h s_v(z) = s_v(x \circ^2_h z)$,~~~
$t_v(x) \circ^2_h t_v(z) = t_v(x \circ^2_h z)$,
\item
$s_h(x) \circ^2_v s_h(y) = s_h(x \circ^2_v y)$,~~~
$t_h(x) \circ^2_v t_h(y) = t_h(x \circ^2_v y)$,
\end{itemize}
for all $(x,z) \in \D \times_{\circ^2_h} \D$ and all $(x,y) \in \D \times_{\circ^2_v} \D$.
 
\medskip

The theorem below shows that in the pure algebraic language, a double category can
be formulated
concisely as a set (of 2-morphisms) with two commuting partial monoid structures.  Note that although
Ehresmann's definition already implies that the structure maps are functors,
the following is not a straightforward translation, especially as the identity assigning map
$\id_{\C} : \Ob(\C) \to \C$ has no explicit presence here. 

\begin{theorem}
A double category (cubical set class one), can be given an equivalent description
as a set $\D$ with two
partial monoid structures $\D_1 = (\D, \circ_1, s_1, t_1)$ and $\D_2 = (\D, \circ_2, s_2, t_2)$
such that the multiplication $\circ_i : \D_j \times \D_j \to \D_j$, $i \neq j \in \{ 1,2 \}$ and the
source and target partial identity maps $s_i, t_i : \D_j \to \D_j$, $i \neq j \in \{ 1,2 \}$ are
all partial monoid homomorphisms. 
\end{theorem}

The non-trivial part of the proof is to unpack the construction to obtain the required four categorical
structures together with all the appropriate compatibility conditions as recalled above. In other
words, we only need to demonstrate that we can construct a double category from the data given in
the statement because it will then be clear that any double category will fit this new
description. See also the diagrams (a),(b) and (c) below.

\begin{proof}

Let $\D$ be a set equipped with two partial monoid structures,\\ $\D_h = (\D, \circ_h, s_h, t_h)$
and $\D_v = (\D, \circ_v, s_v, t_v)$.

From $\circ_h : \D_v \times \D_v \to \D_v$ and equation (1) we find that the two
composition laws $\circ_h$ and $\circ_v$ are compatible with one another but it still
remains to derive the conditions that select the cubical set class before we can unpack the exchange law
below in equation (11) because the cubical set class determines the arrangement of the 2-arrows
when forming their compositions.

Then plugging the multiplication maps $\circ_h$ and $\circ_v$ into equations (2) give compatibility
conditions:
\begin{eqnarray} \label{compatibility}
 s_v(x) \circ_h s_v(z) = s_v(x \circ_h z), ~~~~  t_v(x) \circ_h t_v(z) = t_v(x \circ_h z) \\
 s_h(x) \circ_v s_h(y) = s_h(x \circ_v y), ~~~~ t_h(x) \circ_v t_h(y) = t_h(x \circ_v y)
\end{eqnarray}
 for all $(x,z) \in \D \times_{\circ_h} \D$ and all $(x,y) \in \D \times_{\circ_v} \D$. 

Secondly, the partial identities $s_v, t_v : \D_h \to \D_h$,~~~$s_h, t_h : \D_v \to \D_v$ satisfying
equations (2)  involves
\begin{gather}
\F(s_1(x)) = s_2(\F(x)) ~~~~~~~ => ~~~~~ s_v ( s_h) = s_h (s_v)
\end{gather}
and altogether we have, $\forall x \in \D$, 
\begin{gather}  \label{source}
  s_h (s_v)(x) = s_v (s_h)(x),~~~~
  t_h (t_v)(x) = t_v (t_h)(x)\\
   t_h (s_v)(x) = s_v (t_h)(x), ~~~~
  s_h (t_v)(x) = t_v (s_h)(x).
\end{gather}
and then from the partial identity maps and equation (1), the same set of compatibilities as before
\eqref{compatibility}, (4) are repeated. 

The previous equations establish that the cubical set class is class one, and this means that the
compositions are understood (now that we know which cubical set class the 2-arrows belong to, we can
represent elements using the square diagrams below and hence form their compositions) and we may
now plug $\circ_h : \D_v \times \D_v \to \D_v$ into equation (1) and retrieve the exchange law:
\begin{gather}
(x,y), (w,z) \in \D_v \times \D_v \\
 \circ_h((x,y)(w,z)) = (x \circ_v w) \circ_h (y \circ_v z) \\
 \circ_h(x,y) ~ \circ_h(w,z) = (x \circ_h y) \circ_v (w \circ_h z) \\
 (x \circ_v w) \circ_h (y \circ_v z) = (x \circ_h y) \circ_v (w \circ_h z)
\end{gather}
(whenever both sides of (11) exist).\\
Next we show the emergence of the subsets $\D^1_v$, $\D^1_h$, $\D^0 \subset \D$.\\
The images of the partial identity maps define two distinct subsets
$\D^1_v$ and $\D^1_h$ of 2-arrows, which due to (6),(7) are partial monoids
in their own right, 
\begin{gather}
 \mathrm{im}(t_h)= \mathrm{im}(s_h)~~ =>~~ (\D^1_v, \circ_v, s_v,t_v) \\
 \mathrm{im}(s_v)=\mathrm{im}(t_v)~~ => ~~ (\D^1_h, \circ_h, s_h, t_h)
\end{gather}
By iterating the application of the partial identity maps, one additonally finds partial monoid
homomorphisms:
\begin{gather}
s_v :  (\D^1_v, \circ_v, s_v,t_v) \to (\D^1_v, \circ_v, s_v,t_v), ~~~~
t_v :  (\D^1_v, \circ_v, s_v,t_v) \to (\D^1_v, \circ_v, s_v,t_v), \\
s_h :  (\D^1_h, \circ_h, s_h, t_h) \to (\D^1_h, \circ_h, s_h, t_h), ~~~~
t_h :  (\D^1_h, \circ_h, s_h, t_h) \to (\D^1_h, \circ_h, s_h, t_h),
\end{gather}
whose image sets coincide by equations \eqref{source}, (7). Since all partial identity maps
commute and are idempotent, this iterative procedure can only be performed twice before reaching the
situation
where a subset $\D^0 \subset \D$ is constructed as $\D^0 = \{ x \in \D ~ : ~ s_h(x) = t_h(x) = s_v(x) =
t_v(x) \}$. For clarity and completeness, we state that the elements of $\D^0$ give the identity
morphisms,
or category units in the double category we constructed and that these are compatible with both
partial multiplication rules. Now we see that $\D$ is a double category with $\Ob(\D) =
\D^0$.  
\end{proof}

We can conclude that typical elements of $\D$ are represented by (a) a square 2-cell diagram from
cubical set class one,
\begin{equation*}
(a)
\xymatrix@C30pt@R30pt{
A \ar[r]^{a} \ar[d]_{c} \ar@{.>}[dr]|{x} & B \ar[d]^{b}
\\
C \ar[r]_{d} & D
}
\quad
(b)
\xymatrix@C30pt@R30pt{
A \ar[r]^{a} \ar[d]_{\id^1_v(A)} \ar@{.>}[dr]|{s_v(x)} & B \ar[d]^{\id^1_v(B)}
\\
A \ar[r]_{a} & B
}
\quad
(c)
\xymatrix@C30pt@R30pt{
A \ar[r]^{\id^1_h(A)} \ar[d]_{\id^1_v(A)~} \ar@{.>}[dr]|{s_h(s_v)(x)} & A \ar[d]^{~\id^1_v(A)}
\\
A \ar[r]_{\id^1_h(A)} & A
}
\end{equation*}
(b) a typical element of $\D^1_h$, and (c) by iteration, an object $A \in \D^0$. The diagrams also
illustrate how the data of the maps $\id_v^1 : \Ob(\D) \to \D^1_v$ and $\id_h^1 : \Ob(\D) \to \D^1_h$
(or $\iota : \C^0 \to \C$) have been absorbed into the source and target maps.

\subsection{Other higher structures}

Observe that since the procedure used in the theorem proof is iterative, if $\C$ is a
set with $n$ partial monoid structures $(\C, \circ_i, s_i,t_i)_{i = 1..n}$, we may apply the
same procedure
$n$ times, until $s_i(x)=t_j(x)$ for all $i,j \leq n$.
We propose therefore \emph{a purely algebraic formulation of the notion of (strict) $n$-fold
category}, as a set $\C$, with $n$ partial monoid structures $(\C, \circ_i, s_i, t_i)$, for
$i \in I = \{1..n \}$, in which all
multiplication $\circ_i$ and partial identity maps $s_i$, $t_i$ are partial monoid homomorphisms with
respect to each $j \in I$ such that $j \neq i$. (See \cite{BH} for the definition of an $n$-fold
category.)

\medskip

The following shows that a double category is a
more naturally occurring algebraic object than a 2-category.
A 2-category can be given as a double category in which all horizontal
units $\iota_h$ are also vertical units $\iota_v$ \cite{RR}. Thus a 2-category can be equivalently
described as a set $\C$ equipped with two partial monoid structures $(\C, \circ_1, s_1, t_1)$ and
$(\C, \circ_2, s_2,t_2)$ such that all structure maps ($s_i,t_i,\circ_i$) are partial monoid
homomorphisms and in which $s_1=s_1(s_2)=s_1(t_2)$ and $t_1= t_1(t_2) = t_1(s_2)$ and
$s_1 : \C \to \C^0$ and $t_1: \C \to \C^0$ are surjective maps. 

\medskip

A double category $\D$ with only one object $A \in \D^0$ is equivalent to a set with two
commuting partial monoid structures (in the sense studied above) in which all compositions of the
partial identity maps are constant maps. Thus, $s_h (s_v)(x) = s_v (s_h)(x)$  $= t_h (t_v)(x) = t_v
(t_h)(x)$ $= t_h (s_v)(x) = s_v (t_h)(x)$ $= s_h (t_v)(x) = t_v (s_h)(x) = A$ for all $x \in \D$.

\section{Applications and discussion}

Longo and Roberts \cite{LR} introduced the notion of 2-C*-category in terms of an additional structure
on the set of 2-arrows in a
2-category. In particular, their definition of 2-C*-category begins with an underlying set with two
compositions, each making it into a category. A double category approach to these and other involutive
categories in this language of
partial monoids is further motivated in a coming preprint \cite{idc}. In the simplest terms, an
involutive double category can be thought of as a set with two commuting partial
monoid structures $v$ and $h$ (in the precise sense studied above), each equipped with an involution
$*_v$ and $*_h$ such that $*_h*_v=*_v*_h$.

\medskip

The following describes a motivating example of an involutive double category. Let $(\E,
\pi, \G)$ be a complex line bundle over a discrete double pair groupoid $\G$ (that is,
$\G$ is a double groupoid in which each of the four underlying categorical structures is a discrete pair
groupoid). Modulo a choice of strictification,
the set of fibres of $\E$ has the structure of a double
groupoid. Consider now a section $\sigma$ of $\pi$ over a 2-arrow $g \in \G$ in the sense that its
components consist of an element of each of the 9 fibres over the 2-arrows defined by $g$ (that is, $g$,
$s_h(g)$, $s_v (s_h) (g)$ and so on). We hint at a potential connection with the non-commutative
standard model \cite{sap} as follows. Let $u_L \in A$, $d_L \in C$, $u_R \in B$, $d_R \in D$ labelling
components of $\sigma$ in terms of diagram (a), and let $M$ be the restriction of $\sigma$ to the
groupoid $\G_v$. Note that $M + M^*$ provides an automorphism in a category of Fell bundles over
groupoids \cite{BCL En}, (for Fell bundles see \cite{fbg}) and a self-adjoint linear operator $\Df = M +
M^*$ on the Hilbert space $\H=\mathbb{C}^2 \oplus \mathbb{C}^2$. 

\medskip

A localic groupoid is defined in \cite{Pedro} as a groupoid internal in Loc (recall that a topological
groupoid is a groupoid internal in Top). Examples of localic groupoids tend to be constructed from
a set of arrows with the structure of a groupoid (that is, a partial monoid with all inverses). One
might suggest a definition for a localic double groupoid $\L$ as a set with two commuting partial monoid
structures (again, in the precise sense studied above) such that each partial monoid contained in $\L$
has the structure of a locale. 

\medskip

We add the suggestion that the symmetries of a Penrose
tiling can be described by a partial action $\alpha$ of the group $\G = SO(2) \ltimes \mathbb {R}$ on
the underlying inverse category defined by the patterns in the tiling (this inverse category is detailed
in \cite{Tiling smgs}). Since $\G$ is a subgroup of the Poincar\'e group, which Majard \cite{Majard}
describes as a double group, it would be interesting to investigate the structure of set of ultrafilters
of the inverse semigroup arising from the partial action $\alpha$.

\section{Acknowledgements}

Thank you, grazie and obrigada to Paolo Bertozzini, Pedro Resende and Roberto Conti for discussions and
especially to Paolo Bertozzini for the very significant contributions to this paper.

This work was supported by Funda\c{c}$\tilde{\mathrm{a}}$o   para
as Ci\^encias e a Tecnologia (FCT) through the following projects:\\ PEst-OE/EEI/LA0009/2013,
SFRH/BPD/32331/2006, POCI 2010/FEDER (CAMGSD),\\ EXCL/MAT-GEO/0222/2012 (IST Geometry and
mathematical physics project) and PEst-OE/MAT/UIO143/2014 (CAUL).

\end{document}